\documentclass[11pt,reqno]{amsart}

\usepackage{amsfonts,amssymb,latexsym,mathrsfs}

\large\normalsize
\newtheorem{theorem}{Theorem}

\newcommand{\R}{\ensuremath{\mathbb{R}}}

\newcommand{\T}{\ensuremath{\mathbb{T}}}

\numberwithin{equation}{section}
\numberwithin{theorem}{section}

\begin{document}

\title[first-order nonlinear problem]{existence of three solutions for a first-order problem with nonlinear nonlocal boundary conditions}

\author[d. r. anderson]{Douglas R. Anderson}
\address{Department of Mathematics, Concordia College, Moorhead, MN 56562 USA}
\email{andersod@cord.edu}


\begin{abstract} 
Conditions for the existence of at least three positive solutions to the nonlinear first-order problem with a nonlinear nonlocal boundary condition given by
\begin{eqnarray*}
  && y'(t) - p(t)y(t) = \sum_{i=1}^m f_i\big(t,y(t)\big), \quad t\in[0,1], \\
  && \lambda y(0) = y(1) + \sum_{j=1}^n \Phi_j(\tau_j,y(\tau_j)), \quad \tau_j\in[0,1], 
\end{eqnarray*}
are discussed, for sufficiently large $\lambda>1$. The Leggett-Williams fixed point theorem is utilized.
\end{abstract}

\maketitle
\thispagestyle{empty}

\section{Introduction}

We are interested in the first-order boundary value problem with nonlinear nonlocal boundary condition given by
\begin{eqnarray}
  && y'(t) - p(t)y(t) = \sum_{i=1}^m f_i\big(t,y(t)\big), \quad t\in[0,1], \label{eq} \\
  && \lambda y(0) = y(1) + \sum_{j=1}^n \Phi_j(\tau_j,y(\tau_j)), \quad \tau_j\in[0,1], \label{bc}
\end{eqnarray}
where $p:[0,1]\rightarrow[0,\infty)$ is continuous;
the nonlocal points satisfy $0\le\tau_1<\tau_2<\cdots<\tau_n\le 1$;
the nonlinear functions $\Phi_j:[0,1]\times[0,\infty)\rightarrow[0,\infty)$ satisfy
\begin{equation}\label{phij}
 0 \le y\phi_j(t,y) \le \Phi_j(t,y) \le y\psi_j(t,y), \quad t\in[0,1], \quad y\in[0,\infty),
\end{equation}
for some positive continuous (possibly nonlinear) functions $\phi_j,\psi_j:[0,1]\times[0,\infty)\rightarrow[0,\infty)$;
the scalar $\lambda$ satisfies
\begin{equation}\label{betaj}
  \lambda > \left(1+\sum_{j=1}^n \beta_j\right)\exp\left(\int_0^1p(\eta)d\eta\right) > 1, \quad \beta_j:=\max_{[0,1]\times[0,C]} \psi_j(t,y)
\end{equation}
for some real constant $C>0$;
and the nonlinear functions $f_i:[0,1]\times[0,\infty)\rightarrow[0,\infty)$ are all continuous.
We also set
\begin{equation}\label{alfaj}
 \alpha_j:=\min_{[0,1]\times[B,\lambda B]} \phi_j(t,y)
\end{equation}
for some constant real $B>0$ for later reference. Note that by continuity and compactness $\alpha_j$ and $\beta_j$ exist and satisfy $\beta_j>\alpha_j>0$.

Some of the motivation for this paper and the study of problem \eqref{eq}, \eqref{bc} is as follows. First-order equations with various boundary conditions, including multi-point and nonlocal conditions, are of recent interest. For example, we cite the following papers.
Zhao and Sun \cite{zs} were concerned with the first-order PBVP (if $\T=\R$)
\begin{eqnarray}
  && y'(t) + p(t)y(t) = \lambda f\big(t,y(t)\big), \quad t\in[0,1],  \\
  && y(0) = y(1).
\end{eqnarray}
Tian and Ge \cite{tg} investigated the first-order three-point problem (if $\T=\R$)
\begin{eqnarray}
  && y'(t) + p(t)y(t) = \lambda f\big(t,y(t)\big), \quad t\in[0,1],  \\
  && y(0) -\alpha y(\eta) = \gamma y(1),
\end{eqnarray}
while Gao and Luo \cite{gao} were interested in the problem (if $\T=\R$)
\begin{eqnarray}
  && y'(t) + p(t)y(t) = \lambda f\big(t,y(t)\big), \quad t\in[0,1],  \\
  && y(0) = \sum_{j=1}^{n} \gamma_j y(t_j);
\end{eqnarray}
similarly Anderson \cite{a} studied the first-order problem (if $\T=\R$)
\begin{eqnarray}
  && y'(t) + p(t)y(t) = \lambda f\big(t,y(t)\big), \quad t\in[0,1],  \\
  && y(0) = y(1) + \sum_{j=1}^{n} \gamma_j y(t_j).
\end{eqnarray}
In a related paper, Nieto and R. Rodríguez-L\'{o}pez \cite{nr} considered
\begin{eqnarray}
  && y'(t) + p(t)y(t) = \lambda f\big(t\big), \quad t\in[0,1],  \\
  && \lambda y(t_0) = \sum_{j=1}^{n} \gamma_j y(t_j).
\end{eqnarray}
Gilbert \cite{gil} looked at (if $\T=\R$)
\begin{eqnarray}
  && y'(t) = \lambda f\big(t,y(t)\big), \quad a.e.\; t\in[0,1],  \\
  && y(0) = y(1), \quad\text{or}\quad y(0)=y_0 
\end{eqnarray}
using measure theory and $\Delta$-Carath\'{e}odory functions.
Goodrich \cite{go1} analyzed the $p$-Laplacian problem (if $\T=\R$)
\begin{eqnarray}
  && \phi_p\big(y'(t)\big)=h(t)f\left(y(t)\right), \quad t\in[0,1],  \\
  && y(0) = \psi(y) \quad\text{or}\quad y(0) = B_0\big(y'(1)\big) \quad\text{or}\quad y(0) = \big(y'(1)\big)^m,
\end{eqnarray}
while Graef and Kong \cite{gk} explored the related $p$-Laplacian problem (if $\T=\R$)
\begin{eqnarray}
  && \phi_p\big(y'(t)\big)=f\left(t,y(t)\right), \quad t\in[0,1],  \\
  && y(0) = B\big(y'(1)\big).
\end{eqnarray}
Otero-Espinar and Vivero \cite{ov} gave a general view of different kinds of weak first-order discontinuous boundary value
problems on an arbitrary time scale with nonlinear functional boundary value conditions.
Shu and Deng \cite{sd} proved the existence of three positive solutions to (if $\T=\R$)
\begin{eqnarray}
  && y'(t) = \lambda f\big(y(t)\big), \quad t\in[0,1],  \\
  && y(0) = \gamma y(1).
\end{eqnarray}
Zhao \cite{z} applied a monotone iteration method to the problem (if $\T=\R$)
\begin{eqnarray}
  && y'(t) + p(t)y(t) = \lambda f\big(t,y(t)\big), \quad t\in[0,1],  \\
  && y(0) = g(y(1)),
\end{eqnarray}
where $g$ denotes a nonlinear boundary condition.
Precup and Trif \cite{pt} dealt with the existence, localization and multiplicity of positive solutions to non-local problems for first order differential systems of the form
\begin{eqnarray}
  && y'(t) = f\big(t,y(t)\big), \quad t\in[0,1],  \\
  && y(0) = \alpha[y],
\end{eqnarray}
where $\alpha$ is linear and continuous.
In a second-order problem, Goodrich \cite{go2} was concerned with (if $\T=\R$)
\begin{eqnarray}
  && y^{\Delta\Delta}=-\lambda f\left(t,y(t)\right), \quad t\in[0,1],  \\
  && y(0) = \phi(y), \quad y(1)=0.
\end{eqnarray}
where $\phi$ is a nonlocal boundary condition.
One can see from these representative works that problem \eqref{eq} with the nonlinear nonlocal boundary condition \eqref{bc} is new. 


\section{at least three positive solutions}

In this section we establish the existence of at least three positive solutions for the boundary value problem \eqref{eq}, \eqref{bc}, by finding fixed points for operators on cones in a Banach space; the theorem we will use below is the Leggett-Williams fixed point theorem \cite{le}.

A nonempty closed convex set $P$ contained in a real Banach space $S$ is called a cone if it satisfies the following two conditions:
\begin{enumerate}
  \item[$(i)$] if $y\in P$ and $\zeta \ge 0$ then $\zeta y\in P$;
  \item[$(ii)$] if $y\in P$ and $-y\in P$ then $y=0$.
\end{enumerate}
The cone $P$ induces an ordering $\leq$ on $S$ by $x \le y$ if and only if $y-x\in P$. 
An operator $K$ is said to be completely continuous if it is continuous and compact (maps bounded sets into relatively compact sets). 
A map $\theta$ is a nonnegative continuous concave functional on $P$ if it satisfies the following conditions:
\begin{enumerate}
  \item[$(i)$]  $\theta :P\rightarrow [0,\infty )$ is continuous;
  \item[$(ii)$] $\theta (tx+(1-t)y)\geq t\theta (x)+(1-t)\theta (y)$ for all 
                $x,y\in P$ and $0\leq t\leq 1$.
\end{enumerate}
Let 
\begin{equation}\label{pc}
 P_C := \{y\in P: \| y \| \ < C \}
\end{equation}
and 
\begin{equation}\label{ptheta}
 P(\theta,A,B) := \{ y\in P : A\le\theta (y),\; \| y \| \le B \}. 
\end{equation}


\begin{theorem}[Leggett-Williams]\label{LWTHM} 
Let $P$ be a cone in the real Banach space $S$, $K:\overline{P_C}\rightarrow \overline{P_C}$ be completely continuous, and $\theta$ be a nonnegative continuous concave functional on $P$ with $\theta(y)\le\|y\|$ for all $y\in \overline{P_C}$. Suppose there exist constants $0<A<B<B^\dagger\le C$ such
that the following conditions hold:
\begin{enumerate}
 \item[$(i)$] $\{y\in P(\theta,B,B^\dagger): \theta(y)>B\}\neq\emptyset$ and $\theta(Ky)>B$ for all $y\in P(\theta,B,B^\dagger)$;
 \item[$(ii)$] $\|Ky\|<A$ for $\|y\|\le A$;
 \item[$(iii)$] $\theta(Ky)>B$ for $y\in P(\theta,B,C)$ with $\|Ky\|>B^\dagger$.
\end{enumerate}
Then $K$ has at least three fixed points $y_1$, $y_2$, and $y_3$ in $\overline{P_C}$ satisfying: 
$$ \|y_1\|<A, \quad \theta(y_2)>B, \quad A<\|y_3\|\;\; \mbox{with}\;\; \theta(y_3)<B.$$
\end{theorem}

In the next theorem we use the following notation. Let
\begin{equation}\label{green}
 G(t,s)=\frac{\exp\left(\int_s^t p(\eta)d\eta\right)}{\lambda-\exp\left(\int_0^1 p(\eta)d\eta\right)} \times \begin{cases} 
         \lambda & : 0\le s < t \le 1, \\ 
         \exp\left(\int_0^1 p(\eta)d\eta\right) & : 0 \le t \le s \le 1, 
        \end{cases}
\end{equation}
be the corresponding Green's function,
\begin{equation}\label{mdef}
 M:=\frac{1}{\int_0^1 G(1,s)ds}\left(1-\frac{\exp\left(\int_0^1 p(\eta)d\eta\right)\sum_{j=1}^n\beta_j}{\lambda-\exp\left(\int_0^1 p(\eta)d\eta\right)}\right) > 0
\end{equation}
for $\beta_j$ from \eqref{betaj}, and
\begin{equation}\label{ndef}
 N:=\frac{1}{\int_0^1G(0,s)ds}\left(1-\frac{\sum_{j=1}^n\alpha_j}{\lambda-\exp\left(\int_0^1 p(\eta)d\eta\right)}\right) > 0
\end{equation}
for $\alpha_j$ from \eqref{alfaj}. Both $M$ and $N$ are positive by \eqref{betaj}.


\begin{theorem}
Suppose that \eqref{phij} and \eqref{betaj} hold, and suppose that there exist constants $0<A<B<\lambda B \le C$ such that
\begin{enumerate}
 \item[$(F_1)$] $f_i(t,y) < MA/m$ for $t\in[0,1]$, $y\in[0,A]$,
 \item[$(F_2)$] $f_i(t,y) > NB/m$ for $t\in[0,1]$, $y\in[B,\lambda B]$,
 \item[$(F_3)$] $f_i(t,y)\le MC/m$ for $t\in[0,1]$, $y\in[0,C]$,
\end{enumerate}
for $i=1,2,\cdots,m$, where $M$ and $N$ are given above by \eqref{mdef} and \eqref{ndef}, respectively. 
Then the boundary value problem \eqref{eq}, \eqref{bc} has at least three positive solutions $y_1$, $y_2$, $y_3$ satisfying
$$ \|y_1\|=y_1(1)<A, \quad B<\theta(y_2)=y_2(0), \quad \|y_3\|=y_3(1)>A\;\; \mbox{with} \;\; \theta(y_3)=y_3(0)<B, $$
where $\theta$ is given below in \eqref{line8}.
\end{theorem}

\section{closing comments}

Although we deal with problem \eqref{eq}, \eqref{bc} on the real unit interval, the boundary value problem and accompanying techniques introduced in this work can be readily extended to related difference equations and dynamic equations on time scales.

The main conditions here are \eqref{phij} on the nonlinear functions, and on $\lambda$ in \eqref{betaj}. Note that these assumptions are fairly mild. We do not assume that the $\Phi_j$ are completely separable, nor that they are asymptotically linear. As for $\lambda$, we merely have it bounded below, leaving an unbounded range of possible values that it may assume.

Future authors may wish to explore \eqref{eq}, \eqref{bc} using monotone iteration methods, upper and lower solution techniques, existence of solutions under singularity, positone, or semi-positone conditions, and so on. We leave the details to interested readers.


\end{document}